\documentclass[11pt,a4paper,reqno]{amsart}

\usepackage{amsmath}
\usepackage{amsfonts}
\usepackage{a4wide}
\usepackage{amssymb}
\usepackage{amsthm}

\theoremstyle{plain}

\newtheorem{teo}{Theorem}[section]
\newtheorem{lemma}[teo]{Lemma}
\newtheorem{prop}[teo]{Proposition}
\newtheorem{cor}[teo]{Corollary}
\newtheorem{ackn}{Acknowledgments\!}

\theoremstyle{definition}

\theoremstyle{remark}

\newtheorem{rem}[teo]{Remark}

\numberwithin{equation}{section}

\def\SS{{{\mathbb S}}}

\def\RR{{\mathbb R}}

\def\RRR{{\mathrm R}}
\def\WWW{{\mathrm W}}

\def\BBB{{\mathrm B}}

\def\Ric{{\mathrm {Ric}}}

\def\CCC{{\mathrm C}}

\title[A Note on Four dimensional (Anti-)\,Self--Dual Quasi--Einstein Manifolds]{A Note on Four dimensional (Anti--)\,Self--Dual \\Quasi--Einstein Manifolds}
\date{\today}

\author{Giovanni Catino}
\address[Giovanni Catino]{SISSA --- International School for Advanced Studies, \mbox{Via Bonomea 265,   34136 Trieste (IT)}}
\email{catino@sissa.it}

\date{\today}

\begin{document}

\begin{abstract} In this short note we prove that any complete four dimensional anti--self--dual (or self--dual) quasi--Einstein manifold is either Einstein or locally conformally flat. This generalizes a recent result of X. Chen and Y. Wang. 
\end{abstract}

\maketitle

\begin{center}


\medskip

\centerline{\bf AMS subject classification:  53C24, 53C25}

\end{center}

\section{Introduction}

In this note we will generalize a recent result of X.~Chen and Y.~Wang~\cite{ChenWang} concerning four dimensional (anti--)self--dual gradient Ricci solitons to the case of (anti--)self--dual quasi--Einstein manifolds. We recall that a Riemannian manifold $(M^n,g)$, $n\geq 3$, is a {\em quasi--Einstein manifold\/} if there exist a smooth function $f: M^n\to\mathbb{R}$ and  two constants $\mu,\lambda\in\RR$ such that
\begin{equation}\label{qem}
\Ric+\nabla^{2} f - \mu\, df \otimes df = \lambda g \,.
\end{equation}

When $\mu=0$, quasi--Einstein manifolds correspond to gradient Ricci solitons and when $f$ is constant~\eqref{qem} gives the Einstein equation and we call the quasi--Einstein metric trivial. We also notice that, for $\mu=\tfrac{1}{2-n}$, the metric
$\widetilde{g}=e^{-\frac{2}{n-2}f}g$ is Einstein. Indeed, from the expression of the Ricci tensor of a conformal metric, we get
\begin{eqnarray*}
\Ric_{\widetilde{g}}&=&\Ric_{g} +\nabla^{2} f +\tfrac{1}{n-2} df \otimes
df +\tfrac{1}{n-2} \big(\Delta f - |\nabla f|^{2} \big) g\\
&=& \, \tfrac{1}{n-2}\big(\Delta f - |\nabla f|^{2}+(n-2)\lambda\big)
e^{\frac{2}{n-2}f} \,\widetilde{g}\,.
\end{eqnarray*}
Quasi--Einstein manifolds have been recently introduced by J.~Case,
Y.-S.~Shu and G.~Wei in~\cite{caseshuwei}. In that work the authors focus mainly on the case $\mu\geq 0$. The case $\mu=\frac{1}{m}$ for some $m\in\mathbb{N}$, $m\geq 1$, is particularly relevant due to the link with Einstein warped products. Indeed in~\cite{caseshuwei}, following the results in~\cite{kimkim}, such quasi--Einstein metrics can be characterized as base metrics of Einstein warped product metrics. This characterization on the one hand allows one to translate results from one setting to the other, and on the other hand allows one to exhibit several examples of quasi--Einstein manifolds 
(see~\cite[Chapter 9]{Besse},~\cite{lupagepope}). 

As a generalization of Einstein manifolds, quasi--Einstein manifolds exhibit a certain rigidity. This is well known for $\mu=0$, but we have evidence of this also in the case $\mu\neq 0$. For instance, the author with C.~Mantegazza, L.~Mazzieri and M.~Rimoldi~\cite{mancatmazz1} proved that any locally conformally flat quasi--Einstein manifold of dimension $n\geq 3$ is locally a warped product with $(n-1)$--dimensional fibers of constant curvature (see also~\cite{HePetWylie}). 

In this short note we will prove that in dimension four, to have such a local characterization, it is sufficient to assume the quasi--Einstein metric to be half conformally flat. We recall that a metric is half conformally flat if it is self--dual or anti--self--dual, namely if $\WWW^{-}=0$ or $\WWW^{+}=0$, respectively (see~\cite[Chapter 13, Section C]{Besse} for a nice overview on half conformally flat manifolds).

As we have already seen, the case $\mu=-1/2$ is very special, since it implies that the metric is globally conformally Einstein. In particular, if $(M^{4},g)$ is half conformally flat, then $(M^{4}, e^{-f}\,g)$ is half conformally flat and Einstein.

\medskip

Our main result is the following

\begin{teo}\label{teo} Any complete four dimensional half conformally flat quasi--Einstein manifold with $\mu\neq -1/2$ is either Einstein or locally conformally flat. 
\end{teo}

We want to point out that, following the paper by X.~Chen and Y.~Wang~\cite{ChenWang}, to prove Theorem~\ref{teo} we will make use in a crucial way of the techniques introduced by H.-D. Cao and Q. Chen in~\cite{caochen}.

Using the result in~\cite{mancatmazz1}, we obtain

\begin{cor}\label{cor} Let $(M^{4},g)$ be a complete, simply connected, half conformally flat quasi--Einstein manifold. Then $g$ is either Einstein or of the form 
$$
g \, = \, dt^{2} \, + \,h^{2}(t)\, \sigma^{K} \,,
$$ 
where $\sigma^{K}$ is a Riemannian metric with constant curvature $K$.
\end{cor}

\begin{rem} Many examples of non-trivial rotationally symmetric quasi--Einstein manifolds were constructed by C.~B\"ohm~\cite{bohm1, bohm2}. It is also proven that, for every $\mu$, there is a unique rotationally symmetric quasi--Einstein metric with $\lambda=0$ on $\RR^{n}$. The case $\mu=0$ correspond to the well known Bryant soliton. In particular, from Corollary~\ref{cor}, we have that, on $\RR^{4}$, any half conformally flat quasi--Einstein metric with $\lambda=0$ is isometric either to a Ricci flat metric or to the rotationally symmetric one constructed by C.~B\"ohm~\cite{bohm2}.
\end{rem}

If we restrict ourself to the case $\mu> 0$ and $\lambda\geq 0$, we can say more about the structure of half conformally flat, quasi--Einstein manifolds. First of all, it was proved in~\cite{kimkim}, that if $(M^{n},g)$ is a compact quasi--Einstein manifold with $\lambda=0$ then $g$ has to be Ricci flat (this is true also in the case $\mu=0$, i.e. if $g$ is a gradient steady Ricci soliton). Hence, as a consequence of a theorem of Hitchin (see~\cite[Theorem 13.30]{Besse}) one has the following
\begin{prop} Let $(M^{4},g)$ be a compact, half conformally flat, quasi--Einstein manifold with $\mu\geq 0$ and $\lambda=0$. Then $(M^{4},g)$ is either flat or its universal covering is isometric to a $K3$ surface with the Calabi--Yau metric.
\end{prop}

On the other hand, if $\mu>0$ and $\lambda>0$, then any quasi--Einstein manifold is compact (see~\cite{qian}) with positive scalar curvature (see~\cite{caseshuwei}). Hence, if one assume the metric to be half conformally flat, from Theorem~\ref{teo}, we obtain that $(M^{4},g)$ is either Einstein or locally conformally flat. In the first case, using again Hitchin's theorem, one gets that $g$ is isometric to $\SS^{4}$ or $\mathbb{CP}^{2}$ with their canonical metrics. On the other hand, if $(M^{4},g)$ is locally conformally flat, then a theorem of Kuiper~\cite{kuiper} implies that its universal covering is globally conformally equivalent to $\SS^{4}$. 

Thus, as an application of Theorem~\ref{teo}, we have proved the following
\begin{cor} Let $(M^{4},g)$ be a complete, half conformally flat, quasi--Einstein manifold with $\mu> 0$ and $\lambda>0$. Then $(M^{4},g)$ is either isometric to $\SS^{4}$ or $\mathbb{CP}^{2}$ with their canonical metrics or its universal covering is globally conformally equivalent to $\SS^{4}$.\end{cor}

\section{Proof of Theorem~\ref{teo}}

The proof of Theorem~\ref{teo} will follow X.~Chen and Y.~Wang~\cite{ChenWang}. To fix the notations, we recall that the Riemann curvature operator of a Riemannian manifold $(M^n,g)$ is defined as in~\cite{gahula} by
$$
\mathrm{Riem}(X,Y)Z=\nabla_{Y}\nabla_{X}Z-\nabla_{X}\nabla_{Y}Z+\nabla_{[X,Y]}Z\,.
$$ 
In a local coordinate system the components of the $(3,1)$--Riemann 
curvature tensor are given by
$\RRR^{d}_{abc}\tfrac{\partial}{\partial
  x^{d}}=\mathrm{Riem}\big(\tfrac{\partial}{\partial
  x^{a}},\tfrac{\partial}{\partial
  x^{b}}\big)\tfrac{\partial}{\partial x^{c}}$ and we denote by
$\RRR_{abcd}=g_{de}\RRR^{e}_{abc}$ its $(4,0)$--version.

\smallskip

{\em In all the paper the Einstein convention of summing over the repeated 
indices will be adopted.}

\smallskip

With this choice, for the sphere $\SS^n$ we have
${\mathrm{Riem}}(v,w,v,w)=\RRR_{abcd}v^aw^bv^cw^d>0$. The Ricci tensor is obtained by the contraction 
$\RRR_{ac}=g^{bd}\RRR_{abcd}$ and $\RRR=g^{ac}\RRR_{ac}$ will 
denote the scalar curvature. The so called Weyl tensor is then 
defined by the following decomposition formula (see~\cite[Chapter~3,
Section~K]{gahula}) in dimension $n\geq 3$,
\begin{equation}\label{Weyl}
\WWW_{abcd}=\,\RRR_{abcd}+\frac{\RRR}{(n-1)(n-2)}(g_{ac}g_{bd}-g_{ad}g_{bc})
- \frac{1}{n-2}(\RRR_{ac}g_{bd}-\RRR_{ad}g_{bc}
+\RRR_{bd}g_{ac}-\RRR_{bc}g_{ad}).
\end{equation}
The Weyl tensor satisfies all the symmetries of the curvature tensor
and all its traces with the metric are zero, 
as it can be easily seen by the above formula.\\
In dimension three $\WWW$ is identically zero for every Riemannian
manifold, it becomes relevant instead when $n\geq 4$, since its vanishing is equivalent to $(M^n,g)$ being {\em locally
  conformally flat}. In dimension $n=3$, on the other hand, locally conformally flatness is equivalent to the vanishing of the Cotton tensor
$$
\CCC_{abc} = \nabla_{c} \RRR_{ab} - \nabla_{b} \RRR_{ac} -
\tfrac{1}{2(n-1)} \big( \nabla_{c} \RRR \, g_{ab} - \nabla_{b} \RRR \,
g_{ac} \big)\,.
$$ 
When $n\geq 4$ note that one can compute, (see \cite{Besse}), that
\begin{equation*}
\,\nabla^d\WWW_{abcd}=-\frac{n-3}{n-2}\CCC_{abc}.
\end{equation*}

We recall the following lemma~\cite[Lemma 2.1]
{mancatmazz1}

\begin{lemma} Let $(M^{n},g)$ be a quasi--Einstein manifold. Then the following identities hold
\begin{eqnarray}
\label{eq1}&\RRR + \Delta f - \mu |\nabla f|^{2} = n\lambda& \\
\label{eq2}&\nabla_{b} \RRR = 2(1-\mu) \RRR_{ab} \nabla^{a} f +2 \mu R
\nabla_{b} f -2(n-1)\lambda\mu \nabla_{b} f& \\
\label{eq3}&\nabla_{c} \RRR_{ab} - \nabla_{b} \RRR_{ac} = -
\RRR_{cbad} \nabla^{d} f +\mu\big( \RRR_{ab} \nabla_{c} f - \RRR_{ac}
\nabla_{b} f\big) - \lambda\mu\big( g_{ab} \nabla_{c} f - g_{ac}
\nabla_{b} f \big)
\end{eqnarray} 
\end{lemma}

{\em  From now on, we will consider the case $\mu\neq\tfrac{1}{2-n}$. }

\medskip

Let $(M^{4},g)$ be a complete half conformally flat quasi--Einstein manifolds. Without loss of generality we can assume to be in the case $\WWW^{+}=0$.

Using equations~\eqref{eq2}~\eqref{eq3} together with the decomposition formula for the Riemann tensor~\eqref{Weyl}, from the definition of the Cotton tensor, we obtain
\begin{eqnarray}\label{formcott}
\CCC_{abc} &=& -\WWW_{abcd} \nabla^{d}f - \frac{1}{2}(1+2\mu)\big(\RRR_{ac}\nabla_{b} f - \RRR_{ab}\nabla_{c}f \big) \\\nonumber
&& -\frac{1}{6}(1+2\mu)\big(\RRR_{bd}\nabla^{d}f g_{ac} - \RRR_{cd}\nabla^{d}f g_{ab} \big) +\frac{\RRR}{6}(1+2\mu) \big( \nabla_{b} f g_{ac} - \nabla_{c} f g_{ab} \big) \,.
\end{eqnarray}
Following the notation in~\cite{caochen} and~\cite{ChenWang}, if we define the tensor $\BBB$ to be
\begin{eqnarray*}
\BBB_{abc} \, = \, -\frac{1}{2}\big(\RRR_{ac}\nabla_{b} f - \RRR_{ab}\nabla_{c}f \big) -\frac{1}{6}\big(\RRR_{bd}\nabla^{d}f g_{ac} - \RRR_{cd}\nabla^{d}f g_{ab} \big) +\frac{\RRR}{6}\big( \nabla_{b} f g_{ac} - \nabla_{c} f g_{ab} \big) \,,
\end{eqnarray*}
which was introduced in~\cite{caochen}, from~\eqref{formcott}, we obtain that
$$
\CCC_{abc} \, = \, - \WWW_{abcd} \nabla^{d} f + (1+2\mu) \, \BBB_{abc} \,.
$$
We recall that we are assuming $\mu\neq -\tfrac{1}{2}$ and we observe that the expression of the tensor $\BBB$ is exactly the same as the one in the proof of Lemma 2.1 in~\cite{ChenWang}. Therefore, following their computation, we can easily prove that, at every point $p\in M^{4}$ where $\nabla f_{|p} \neq 0$, one has $\BBB_{|p}=0$. From this fact it follows that, at $p$, the Ricci tensor either has a unique eigenvalue $\eta_{1}$, or has two distinct eigenvalues $\eta_{1}$ and $\eta_{2}$ of multiplicity 1 and 3 respectively. In either case, $e_{1}=\nabla f / |\nabla f| $ is an eigenvector with eigenvalue $\eta_{1}$. Moreover, as in the proof of Lemma 2.3 in~\cite{ChenWang}, one can easily check that $\WWW_{|p} =0$, whenever $\nabla f _{|p} \neq 0$.

\medskip

Now we can conclude the proof of Theorem~\ref{teo}.

\medskip

Let us assume that the set $\{p\in M^{4}\,|\,|\nabla f|^{2}_{|p} \neq 0\}$ is dense in $M^{4}$. Then, from the previous observations, we have that the Weyl tensor has to vanish on a dense set, which clearly implies that $g$ has to be locally conformally flat.

If the above does not hold, it implies that the function $f$ is constant in an open set of $M$ and $g$ is Einstein in this open set. As it was observed in~\cite[Proposition 2.8]{HePetWylie}, quasi--Einstein metrics, as well as the function $f$, have to be real analytic. Therefore $g$ has to be Einstein on the whole manifold $M^{4}$.

\begin{ackn} The author is partially supported by the Italian project FIRB--IDEAS ``Analysis and Beyond''.
\end{ackn}

\bibliographystyle{amsplain}
\bibliography{AntiselfdualQE}

\bigskip

\parindent=0pt

\end{document}